\documentclass{siamart190516}
\usepackage{lipsum}
\usepackage{amsfonts}
\usepackage{graphicx}
\usepackage{epstopdf}
\usepackage{algorithmic}
\ifpdf
  \DeclareGraphicsExtensions{.eps,.pdf,.png,.jpg}
\else
  \DeclareGraphicsExtensions{.eps}
\fi

\newsiamremark{hypothesis}{Hypothesis}
\crefname{hypothesis}{Hypothesis}{Hypotheses}
\newsiamthm{claim}{Claim}

\headers{Adaptive Two-Layer ReLU Neural Network}{M. Liu, Z. Cai, and J. Chen}

\title{Adaptive Two-Layer ReLU Neural Network:\\
I. Best Least-squares Approximation\thanks{This work was supported in part by the National Science Foundation
under grant DMS-2110571.}}
\author{Min Liu\thanks{School of Mechanical Engineering, Purdue University, 585 Purdue Mall, West Lafayette, IN 47907-2088(\email{liu66@purdue.edu}) }
\and Zhiqiang Cai\thanks{Department of Mathematics, Purdue University, 150 N. University Street, West Lafayette, IN 47907-2067 
  (\email{caiz@purdue.edu} and \email{chen2042@purdue.edu}).}
  \and Jingshuang Chen\footnotemark[3]
  } %\email{chen2042@purdue.edu},
\usepackage{amsopn}

\ifpdf
\hypersetup{
  pdftitle={FOSLS},
  pdfauthor={M.Liu, Z. Cai and J.Chen}
}
\fi

\DeclareUnicodeCharacter{2212}{-}
\usepackage[utf8]{inputenc}
\usepackage{bm}
\usepackage{mathtools}
\usepackage{indentfirst}
\usepackage{subfigure}
\usepackage{tcolorbox}
\usepackage{color}
\usepackage[export]{adjustbox}
\usepackage{makecell}

\usepackage{algorithm }
\usepackage{algorithmic}
\Crefname{ALC@unique}{Line}{Lines}

\newcommand{\R}{\mathbb{R}}

\usepackage{graphicx}
\usepackage{diagbox}
\usepackage{pifont}
\newcommand{\vertiii}[1]{{\left\vert\kern-0.25ex\left\vert\kern-0.25ex\left\vert #1 
    \right\vert\kern-0.25ex\right\vert\kern-0.25ex\right\vert}}

\newcommand{\btheta}{\mbox{\boldmath${\theta}$}}

\newcommand{\bomega}{\mbox{\boldmath${\omega}$}}

\newcommand{\bb}{{\bf b}}
\newcommand{\bv}{{\bf v}}

\newcommand{\bx}{{\bf x}}

\newcommand{\cI}{{\cal I}}
\newcommand{\cK}{{\cal K}}
\newcommand{\cM}{{\cal M}}
\newcommand{\cQ}{{\cal Q}}
\newcommand{\cS}{{\cal S}}
\newcommand{\cT}{{\cal T}}

\begin{document}

\maketitle
\begin{abstract}
In this paper, we introduce adaptive network enhancement (ANE) method for the best least-squares approximation using two-layer ReLU neural networks (NNs). For a given function $f(\bx)$, 
the ANE method generates a two-layer ReLU NN and a numerical integration mesh such that the approximation accuracy is within the prescribed tolerance. The ANE method provides a natural process for obtaining a good initialization which is crucial for training nonlinear optimization problems.
Numerical results for functions of two variables exhibiting either intersecting interface singularities or sharp interior layers demonstrate efficiency of the ANE method.  
\end{abstract}

% REQUIRED

% REQUIRED

\section{Introduction}

Deep neural networks (DNNs) have achieved astonishing performance in computer vision, natural language processing, and many other artificial intelligence tasks. This success encourages wide applications to other fields, including recent studies of using DNN models to numerically solve partial differential equations (PDEs).
Despite their great successes in many practical applications, it is widely accepted that approximation properties of DNNs are not yet well-understood and that understandings on why and how they work could lead to significant improvements. This explains rapidly increasing interests in theoretical and algorithmic studies of DNNs during recent years. 

DNNs produce a new class of functions through compositions of linear transformations and activation functions. Their studies and applications may be traced back to the work of Hebb \cite{Hebb1949} in the late 1940's and Rosenblatt \cite{Rosenblatt1958} in the 1950's. 
An often cited theoretical results on DNNs is the so-called universal approximation property \cite{Cybenko1989, HornikS1989}, e.g., a two-layer NN is dense in $C(\Omega)$ for any compact subset $\Omega\in \R^d$ provided that the activation function is not a polynomial. Moreover, order of approximation for functions in the Sobolve space have been obtained for two-layer NNs using various activation functions \cite{Petrushev1998}. For results on approximation theory of DNNs before 2000, see a survey article by Pinkus \cite{pinkus1999} and references therein.

Despite many efforts and much impressive progress made by numerical analysts, computational scientists, and practitioners, approximation properties of DNNs remain an active and open research field. Without complete understanding of approximation properties of DNNs, current methods on design of network structures are empirical. Tuning of depth and width is tedious, mainly from experimental results in ablation studies which typically require domain knowledge about the underlying problems.
This leads to a fundamental, open question in machine learning: {\it given a target function/PDE, what is the minimal network model required, in terms of width, depth, and number of parameters, to approximate the function/solution within the prescribed accuracy?}

The purpose of this paper is to introduce and study adaptive network enhancement (ANE) methods for the best least-squares approximation to a target function by a two-layer ReLU NN, and, hence, to address this open problem partially.
Specifically, for a given target function $f(\bx)$ and a given tolerance $\epsilon>0$, the ANE method generates a two-layer ReLU neural network such that the approximation accuracy is within the prescribed tolerance. One of key components of the ANE method for the best least-squares approximation to a given function is the enhancement strategy which determines how many new neurons to be added, when the current approximation is not within the given accuracy. To address this issue, we introduce a global and a local network enhancement strategies. The global enhancement is based on a fixed convergence rate (see (\ref{global-n}); and the local one
%The key ingredient of the method is the adaptive network enhancement strategy which determines how many new neurons to be added, when the current approximation is not within the given accuracy. 
is done through local error indicators collected on the physical subdomains plus a proper neuron initialization (detailed in section~5). The ANE method for solving elliptic PDEs is presented in the companion paper \cite{LiuCai2020}.

Another important ingredient is the numerical integration mesh for evaluating the loss function. For many problems in machine learning, integral of the $L^2(\Omega)$ norm is often computed numerically by stochastic sampling approach, which in turn leads to {\it theoretical} convergence rate independent of the dimension. Other numerical integration methods that are independent of the dimension include quasi-Monte Carlo method \cite{DickKuoSloan2014} and the sparse grid method \cite{BungartzGriebel2004}. %, even though its variance of the integrand may be very large. 
For simplicity, 
in this paper, we use adaptive numerical integration based on ``mid-point'' quadrature on either uniform or composite mesh. The composite mesh here means those meshes obtained from adaptive mesh refinement (AMR), where refinement of an element is done by subdividing it into small uniform elements. 
The AMR method presented in the paper is suitable for low dimensional problems and may be replaced by any adaptive integration procedure such as adaptive version of Monte Carlo, quasi-Monte Carlo, or sparse grid, etc. if a high dimensional problem is considered.

Theoretically, we show that the total approximation error is bounded by the approximation error of the NN plus the error of numerical integration (see Theorem~4.1) under the assumption of the Marcinkiewicz problem. This indicates that numerical integration should be chosen to ensure at least the accuracy of the current NN. For simple problems, one may simply use a fine uniform mesh which is able to capture all local behaviors of the integrand. % and then apply the ANE method described in the previous paragraph. 
For computationally intensive problems, one might need to use local AMR to generate a proper composite mesh. The stopping criterion for the AMR is based on if the mesh refinement of numerical integration improves the approximation accuracy (see Algorithm~5.2). With AMR for numerical integration, the ANE method defined in Algorithm~5.3 is able to generate a two-layer ReLU NN and a composite numerical integration mesh such that the approximation accuracy is within the prescribed tolerance. 

The values of the parameters are trained by iteratively ``solving'' the non-convex optimization problem in (\ref{L2App-d}). This high dimensional, non-convex optimization problem tends to be computationally intensive and complicated. Currently, it is often solved by iterative optimization methods such as gradient descent (GD), Stochastic GD, Adam, etc. (see, e.g., \cite{BoCuNo2018} for a review paper in 2018 and references therein). Usually nonlinear optimizations have many solutions, and the desired one is obtained only if we start from a close enough first approximation. The ANE method provides a natural process for obtaining a good initialization. Starting with a relatively small NN, the approximation of the previous NN is already a good approximation to the current NN in the loops of the ANE method. To provide a better approximation than the previous one, we divide all network parameters into two groups: linear parameters (output layer weights and bias) and nonlinear parameters (hidden layer weights and biases). Initialization of nonlinear parameters are based on their physical partitioning of the domain and initial of linear parameters are obtained by solving a system of linear equations with given nonlinear parameters. 

The paper is organized as follows. Section~2 presents two-layer ReLU NNs. The best least-squares approximation and its discrete counterpart are described in sections~3 and 4, respectively. The ANE method is introduced in section 5, and initialization of parameters at different stage are proposed in section~6. Finally, numerical experiments for functions with intersecting interface singularities and interior layer like discontinuities are given in section~7, and conclusion in section~8.

%As seen in the previous section, a two-layer ReLU network produces linear splines with free knots. It is well-known that when the knots of linear splines are free, approximation to functions can be improved dramatically (see, e.g., \cite{jupp1978}). In particular, ``rough'' functions may be approximated without exhibiting overshoots and oscillations.

%Free knot splines have not been as popular as one might expect simply because of the difficulty of numerically solving non-convex minimization problem (\ref{L2App}). Efficient iterative solution methods have been studied by researchers: e.g., Jupp \cite{jupp1978} on a relaxed Gauss-Newton method, Loach-Wathen \cite{wathen1991} initial reference and iterative solvers ...

\section{Two-Layer ReLU Neural Network}

%This paper considers the approximation of functions by two-layer neural networks. 
%In particular, we restrict our attention to two-layer neural networks.
%functions of one variable and two-layer ReLU neural network. Generalizations to functions of several variables and other activation functions may be foreseen. 
A two-layer NN consists of an input and an output layers. The output layer does not have an activation function. Layers other than the output layer are called hidden layers. So a two-layer NN is also referred to as a one-hidden layer NN. 

In $d$-dimension, for $i=1,\,2,\,...,\,n$, let $\bomega_i\in\R^d$ 
and $b_i\in \R$ be the weights and bias of the first (input) layer, respectively; and let $c_i\in\R$ and $c_0\in\R$ be the respective weights and bias of the second (output) layer. Then a two-layer ReLU NN with $n$ neurons produces the following set of functions:
\[
 \hat{\cal M}_n(\sigma) = \left\{c_0+\sum_{i=1}^n c_i\sigma(\bomega_i\cdot \bx -b_i)\, :\,  
c_i,\, b_i\in \R,\,\, \bomega_i\in \R^d \right\},
 \]
where $\sigma$ is the rectified linear unit (ReLU) activation function given by
 \[
 \sigma(t) = \max \{0,\,t\}
 =\left\{\begin{array}{ll}
   0, &  t < 0 ,\\[2mm]
   t , & t \ge 0,
   \end{array}\right.
 \]
for any $t\in\R$. The $\sigma(t)$ is a continuous piece-wise linear function having a breaking point at $t=0$ and belongs to a class of activation functions of the form
\[
\sigma_k(t) = \big(\max \{0,\,t\}\big)^k
 =\left\{\begin{array}{ll}
   0, &  t < 0 ,\\[2mm]
   t^k , & t \ge 0
   \end{array}\right.
   \quad\mbox{for }\,\, k\in \mathbb{Z}_{+},
\]
where $\mathbb{Z}_{+}$ is the set of all positive integers. Note that $\sigma_k(t)\in C^{k-1}(\R)$ is a piece-wise polynomial of degree $k$ with a breaking point at $t=0$. For simplicity of presentation, we restrict our attention to the ReLU activation function. Extension of results in this paper
to general activation functions $\sigma_k(t)$ is straightforward.

There are $(d+2)n+1$ parameters for functions in the set $\hat{\cal M}_n(\sigma)$, where $n+1$ of them are the output weights and bias $\{c_i\}_{i=0}^n$ and $(d+1)n$ of them are the input weights $\{\bomega_i\}_{i=1}^n$ and bias $\{b_i\}_{i=1}^n$. We refer to the former as linear parameters and the later nonlinear parameters.
Thus, $\hat{\cal M}_n(\sigma)$ has of $n+1$ linear
and $(d+1)n$ nonlinear parameters. To remove $n$ nonlinear parameters, we notice that 
\[
 \sigma(\bomega\cdot \bx -b) = |\bomega | \, \sigma\left(\dfrac{\bomega}{|\bomega |}\cdot \bx -\dfrac{b}{|\bomega |}\right),
 \]
where $|\bomega |=\sqrt{\omega_1^2+\cdots + \omega_d^2}$ is the length of a vector $\bomega\in \R^d$. This implies that $\hat{\cal M}_n(\sigma)$ is equal to 
\begin{equation}\label{ReLU-n}
 {\cal M}_n(\sigma,d) = \left\{c_0+\sum_{i=1}^n c_i\sigma(\bomega_i\cdot \bx -b_i)\, :\,  
c_i,\, b_i\in \R,\,\, \bomega_i\in \cS^{d-1} \right\},
 \end{equation}
where $\cS^{d-1}$ is the unit sphere in $\R^d$. The number of parameters in ${\cal M}_n(\sigma,d)$ is \[M(n,d)= (d+1)n +1.\]

%Let $\{w^1_i\}_{i=1}^n$ and $\{b_i\}_{i=1}^n$ be the weights and bias of the first (input) layer, respectively; and let $\{w_i\}_{i=1}^n$ and $w_0$ be the respective weights and bias of the second (output) layer. Then the two-layer ReLU neural network with $n$ neurons produces  functions of the following form
%  \begin{equation}\label{NN}
% \hat{\mathcal{N}}(x %, \btheta
% )=w_0 + \sum_{i=1}^n w_i \sigma(w_i^1 x - b_i) \in \R.
 %\end{equation}
%If $x\in \R^d$, then $w^1_i\in \R^d$ and $b_i\in \R$.By the universal approximation theorem (see, e.g., \cite{pinkus1999}), the class of functions of the form in (\ref{NN}) for all integers $n\in\mathbb{N}$ is dense in the space of continuous functions defined on a compact set $K\subset \R^d$ if and only if $\sigma$ is non-polynomial.

%For a perceptive understanding of ${\cal M}_n(\sigma,d)$, 

Below let us look at ${\cal M}_n(\sigma,d)$ in one-, two- and $d$-dimension, separately. When $d=1$, we have $\cS^{0}=\{-1,\,1\}$. Without loss of generality, we will choose $\omega_i=1$ for all $i=1,\, ...,\,n$. Then %This leads to the following set of functions:
\begin{equation}\label{ReLU-1}
 {\cal M}_n(\sigma,1) = \left\{v(x , {\small\btheta})=c_0+\sum_{i=1}^n c_i\sigma(x -b_i)\, :\,  
c_i,\, b_i\in \R \right\},
 \end{equation}
where ${\small\btheta}=({\bf c}, \, {\bf b})$ denotes all parameters ${\bf c}=(c_0,\,c_1,\, ...,\, c_n)$ and ${\bf b}=(b_1,\, ...,\, b_n)$.
The $\cM_n(\sigma,1)$ is the set of linear splines with $n$ free knots that had been studied intensively in the late 1960s (see, e.g., \cite{Schumaker}). It has been shown that the approximation of functions by linear splines can generally be dramatically improved if the knots are free \cite{jupp1978}; particularly, the Gibbs phenomena for ``rough'' functions can be avoided \cite{Baker85}.

%Since $\sigma(x - b_i)$ %for all $x\in \R$ is a piece-wise linear function having a breaking point at $x=b_i$: 
%\[\sigma(x - b_i)= \max \{0,\,x - b_i\}   =\left\{\begin{array}{ll}   0, & x < b_i ,\\[2mm]
 %  x - b_i, & x \ge b_i   \end{array}\right. \]
%it is then easy to see that ${\cal M}_n(\sigma)$ is the set of continuous piece-wise linear functions with $n$ free breaking points. 
%The total number of parameters of ${\cal M}^1_n(\sigma)$ is $2n+1$ with $n+1$ linear parameters $\{c_i\}_{i=0}^n$ and $n$ nonlinear parameters $\{b_i\}_{i=1}^n$. 
%Denote vectors of weights and bias by
% \[ {\bf c}=(c_0,\,c_1,\, ...,\, c_n)^t \quad\mbox{and}\quad  {\bf b}=(b_1,\, ...,\, b_n)^t, \]
%respectively, then function $v$ in ${\cal M}_n(\sigma)$ can be represented as follows:
%\begin{equation}\label{N} v(x , {\small\btheta})=c_0 + \sum_{i=1}^n c_i \sigma(x - b_i)= \sum_{i=0}^n c_i \varphi_i(x,\, b_i),\end{equation}
%where $\varphi_0(x,\, b_0)=1$, $\varphi_i(x,\, b_i) = \sigma(x-b_i)$ for $i=1,\, ...,\,n$, and ${\small\btheta}=({\bf c}, \, {\bf b})$ are parameters. 
%Since $\mathcal{N}(x , \btheta)$ is linear with respect to the output weights ${\bf c}$ and nonlinear on the input bias ${\bf b}$, the parameters ${\bf c}$ and ${\bf b}$ are referred to linear and nonlinear parameters, respectively. This restricted network has only a half of the number of nonlinear parameters but generates the same approximation class as the standard one.

In two dimensions ($d=2$), $\cS^{1}$ is a unit circle:
 \begin{eqnarray*}
 \cS^{1}
 =\left\{\bomega=(\omega_1,\,\omega_2)^t\in \R^2\, :\, \omega_1^2+\omega_2^2=1\right\} %\\[2mm]
 =\left\{\bomega=\big(\cos\gamma,\,\sin\gamma \big)^t\, :\, 0\leq \gamma\leq 2\pi\right\}.
 \end{eqnarray*}
This gives
\begin{equation}\label{ReLU-2}
 {\cal M}_n(\sigma,2) = \left\{c_0+\sum_{i=1}^n c_i\sigma\big(
 (\cos \gamma_i) \,x_1 +(\sin \gamma_i)\,x_2 -b_i\big) :\, 
c_i,\, b_i\in \R,\,\, \gamma_i\in [0,\,2\pi] \right\},
 \end{equation}
which is the set of continuous piece-wise linear functions with $n$ free lines  
 \begin{equation}\label{lines}
 {l}_i:\,\,(\cos \gamma_i) \,x_1 +(\sin \gamma_i)\,x_2-b_i=0
 \quad\mbox{for } i=1,\,...,\,n.
 \end{equation}
Similarly, in the $d$-dimension, ${\cal M}_n(\sigma,d)$ is the set of continuous piece-wise linear functions with $n$ free hyper-planes
 \begin{equation}\label{planes}
 {\cal P}_i:\, \,\bomega_i\cdot \bx-b_i=0
 \quad\mbox{for } i=1,\,...,\,n.
 \end{equation}
Clearly, ${\cal M}_n(\sigma,d)$ for $d\ge 2$ may be treated as a non-standard but beautiful extension of linear splines with free knots ${\cal M}_n(\sigma,1)$ to multi-dimension.

%Let \[{\cal M}(\sigma) = \mbox{span}\, \left\{\sigma(\bomega\cdot \bx -b)\, :\,  b\in \R,\,\, \bomega\in \R^d \right\}. \]
%The often cited universal approximation property of two-layer NNs \cite{Cybenko1989, HornikS1989, pinkus1999} is that the linear space ${\cal M}(\sigma)$ is dense in $C(K)$, the space of all continuous functions defined on a compact set $K\in \R^d$. 

Let 
\[
%\varphi_0(\bx;\bw_0,\, b_0)=1
%\quad\mbox{and}\quad
\varphi_i(\bx) = \varphi_i(\bx;\bomega_i,\, b_i) = \sigma(\bomega_i\cdot\bx-b_i),
%\quad\mbox{for }\,\, i=1,\, ...,\,n.
\]
%and let $\varphi_0(\bx;\bw_0,\, b_0)=1$.
%Since $\sigma(t)$ is a piece-wise linear function with a breaking point at $t=0$, 
which is a piece-wise linear function with free hyper-planes: $\bomega_i\cdot \bx =b_i$ for $i=1,\, ...,\,n$. Let $\varphi_0(\bx)=\varphi_0(\bx;\bomega_0,\, b_0)=1$. For convenience of readers, we state and prove a well-known fact on the linear independence of $\{\varphi_i(\bx)\}^n_{i=0}$.

\begin{lemma}\label{l-ind}
Assume that hyper-planes $\{\bomega_i\cdot \bx =b_i\}_{i=1}^n$ are distinct. Then $\{\varphi_i(\bx;\bomega_i,\, b_i)\}_{i=0}^n$ are linearly independent.
\end{lemma}
\begin{proof}
Linear independence of $\varphi_0=1$ and $\varphi_1(\bx;\bomega_1,b_1)$ is a direct consequence of the fact that $\varphi_1(\bx;\bomega_1, b_1)\equiv 0$ on $\R^d\setminus \mbox{suppt}\{\varphi_1\}$.
Assume that the lemma is valid for $n=k$, then linear independence of $\{\varphi_i(\bx;\bomega_i, b_i)\}_{i=0}^{k+1}$ follows from the fact that $\sum\limits_{i=0}^{k}c_i \varphi_i(\bx;\bomega_i, b_i)\equiv 0$ for all $\bx \in\R^d\setminus \mbox{suppt}\{\varphi_{k+1}\}$ and the assumption that all hyper-planes $\{\bomega_i\cdot \bx =b_i\}_{i=1}^{k+1}$ are distinct.
This completes the proof of the lemma by induction.
\end{proof}

\section{The Best Least-squares Approximation} % by ReLU Network}

Denote vectors of weights and bias by
 \[
 {\bf c}=(c_0,c_1, ..., c_n),
 \quad \bomega=(\bomega_1, ... , \bomega_n),
 \quad\mbox{and}\quad 
 {\bf b}=(b_1, ..., b_n),
 \]
respectively, then
each function $v\in {\cal M}_n(\sigma,d)$ may be represented as follows:
\begin{equation}\label{N}
 v(\bx ; {\small\btheta})=c_0 + \sum_{i=1}^n c_i\, \sigma(\bomega_i\cdot\bx - b_i)
 = \sum_{i=0}^n c_i \,\varphi_i(\bx; \bomega_i,b_i),
 \end{equation}
where ${\small\btheta} = \big({\bf c},\hat{\small\btheta}\big)$ with $\hat{\small\btheta}=\big(\bomega, \bb\big)$ are parameters. 

For a given function $f(\bx)$ defined on $d$-dimensional domain $\Omega\in\R^d$, the best least-squares approximation is to find $f_n(\bx ; {\small\btheta}^*)\in {\cal M}_n(\sigma,d)$ such that
 \begin{equation}\label{L2App}
  \|f(\cdot)-f_n(\cdot;{\small\btheta}^*)\|=\min_{v\in {\cal M}_n(\sigma,d)} \|f-v\|
  =\min_{{\scriptsize\btheta}\in\R^{M(n,d)}} \|f(\cdot)-v(\cdot;{\small\btheta})\|,
  \end{equation}
where $\|\cdot\|$ denotes the $L^2(\Omega)$ norm, $M(n,d)$ is the number of parameters defined in the previous section, and $v(\bx;{\small\btheta})$ is given in (\ref{N}). It was proven by Petrushev in \cite{Petrushev1998}
(see also \cite{pinkus1999}) that for any $f(\bx)$ in the Sobolev space $H^m(\Omega)$ for $m=1,\,  ... ,\, 2+\dfrac{d-1}{2}$, there exists a positive constant $C$ such that 
 \begin{equation}\label{error}
 \|f-f_n\|\leq C\,n^{-m/d} \|f\|_{H^m(\Omega)}.
 \end{equation}

\begin{remark}
In one dimension, when $f\in L^p(\Omega)$ for $0<p\leq \infty$, it was shown
{\em (}see Rice {\em \cite{Rice1969}} and Powell {\em \cite{Powell1968}}) that problem {\em (\ref{L2App})} has a solution $f_n \in C[0,\,1]$. Solution of problem {\em (\ref{L2App})} is not unique in general; but it is unique for sufficiently smooth $f$ and large enough $n$ {\em (}see Chui et al. {\em \cite{chui1977})}.
\end{remark}

Generally,  ${\cal M}_n(\sigma,d)$ is only a set
of functions. But for a fixed parameter $\hat{\small\btheta}^0=\big(\bomega^0, \bb^0\big)$, the set ${\cal M}_n(\sigma,d)$ becomes a subspace
 \[
 {\cal M}_n(\sigma,d) = \mbox{span}\, \{\varphi_i(\bx; \bomega^0_i,b^0_i)\}_{i=0}^n.
 \]
Then the best least-squares approximation in (\ref{L2App}) becomes to find $f_n^0 = \sum\limits_{i=0}^n c^0_i \,\varphi_i(\bx; \bomega^0_i,b^0_i)\in {\cal M}_n(\sigma,d)$ such that
 \[
 (f_n^0, \varphi_i(\bx; \bomega^0_i,b^0_i)) = (f,\varphi_i(\bx; \bomega^0_i,b^0_i)) \quad\forall\,\, i=0,1, ..., n,
 \]
where $\big(f,\,g\big) =\int_\Omega f(\bx) g(\bx) \,d\bx$ denotes the $L^2(\Omega)$ inner product.
The corresponding system of algebraic equations is 
 \begin{equation}\label{O_W}
 {\bf M} (\hat{\small\btheta}^0)\, {\bf c}^0 = F  (\hat{\small\btheta}^0),
\end{equation}
where ${\bf M}(\hat{\small\btheta}^0) = \big( M_{ij}\big)_{(n+1)\times (n+1)}$ is the mass matrix with $M_{ij}= \big(\varphi_j(\bx; \bomega^0_j,b^0_j), \varphi_i(\bx; \bomega^0_i,b^0_i)\big)$, ${\bf c}^0= (c_0^0,c^0_1, ..., c_n^0)^t$, and $F  (\hat{\small\btheta}^0)= \big(F_i\big)_{(n+1)\times 1}$ is the right-hand side vector with $F_i= \big(f,\varphi_i(\bx; \bomega^0_i,b^0_i)\big) $.

\begin{lemma}
Assume that the hyper-planes $\{\bomega^0_i\cdot\bx=b^0_i\}_{i=1}^n$ are distinct. 
Then the mass matrix ${\bf M} (\hat{\small\btheta}^0)$ is symmetric, and positive definite.
\end{lemma}

\begin{proof}
Clearly, ${\bf M} (\hat{\small\btheta}^0)$ is symmetric.
For any $\bv = (v_0,\,v_1,\, ...,\, v_n)^t$,
we have 
\[
\bv^t {\bf M} (\hat{\small\btheta}^0) \bv = \|v\|^2,
\]
where $v(\bx) = \sum\limits_{i=0}^n v_i\varphi_i(\bx;\,\bomega^0_i,\, b^0_i)$.
By Lemma~2.1, $\|v\|^2$ is positive for any nonzero $\bv$, which, in turn, implies that ${\bf M} (\hat{\small\btheta}^0)$ is positive definite.
\end{proof}

\section{Effect of Numerical Integration}

In practice, integral of the loss function is often computed numerically. A common practice in machine learning (see, e.g., \cite{Dissanayake94, Sirignano18, Karniadakis19}) uses Monte Carlo integration of the form
\begin{equation}\label{sampling}
 \cI (v) = \int_\Omega v(\bx) \,d\bx 
 \approx \dfrac{|\Omega|}{N} \sum_{i=1}^N v(\bx_i),
 \end{equation}
where $|\Omega|$ is the volume of the domain $\Omega$ and
$\{\bx_i\}_{i=1}^N$ are the sampling points randomly generated based on an assumed distribution of $\bx$. 
This stochastic approach is simple and valid for any dimensions. Moreover, it leads to theoretical convergence rate independent of the dimension. Other numerical integration methods that are independent of the dimension include quasi-Monte Carlo method \cite{DickKuoSloan2014} and the sparse grid method \cite{BungartzGriebel2004}.

In this paper, we use adaptive numerical integration as in \cite{cai2020} in line with the ANE method.  For simplicity of presentation, we consider only ``mid-point'' quadrature on either uniform or composite mesh. The composite mesh here means those meshes obtained from adaptive mesh refinement (AMR), where refinement of an element is done by subdividing it into small uniform elements.  
To this end, let 
\[
{\cal T}=\{K\, :\, K\mbox{ is an open subdomain of } \Omega\}
\]
be a partition of the domain $\Omega$. Here, the partition means that union of all subdomains of ${\cal T}$ equals to the whole domain $\Omega$ and that any two distinct subdomains of ${\cal T}$ have no intersection; more precisely,
 \[
 \bar{\Omega} = \cup_{K\in {\cal T}} \bar{K}
 \quad\mbox{and}\quad
 K\cap T = \emptyset,
 \quad \forall\,\, K,\, T \in {\cal T}.
 \]
Let $\bx_{_T}$ be the centroid of $T\in {\cal T}$. The $\bx_{_T}$ will be used as quadrature points which are fundamentally different from sampling points used in the setting of standard supervised learning. The composite ``mid-point'' quadrature rule is given by
 \[
 \cI (v) = \int_\Omega v(\bx) \,d\bx 
 \approx \sum_{T\in \cT} v(\bx_{_T}) \, |T|
 \equiv \cQ_{_\cT} \big(v\big),
 \]
where $|T|$ is the volume of element $T\in \cT$. Similarly, one may use any quadrature rule such as composite trapezoidal, Simpson, Gaussian, etc.

Let $\cQ_{_\cT}$ be a quadrature operator, i.e., 
 $\cI (v) \approx \cQ_{_\cT} \big(v\big)$, such that 
 \[
 \|v\|_{_\cT}=\sqrt{(v,v)_{_\cT}}=\sqrt{\cQ_{_\cT} \big(v^2\big)}
 \]
defines a weighted $l_2$-norm.
The best discrete least-squares approximation with numerical integration over the partition $\cT$ is to find $f_{_\cT}(\bx ; {{\small\btheta}^*_{_\cT}})\in {\cal M}_n(\sigma,d)$ such that
 \begin{equation}\label{L2App-d}
  \|f(\cdot)-f_{_\cT}(\cdot ; {{\small\btheta}^*_{_\cT}})\|_{_\cT}
  =\min_{v\in {\cal M}_n(\sigma,d)} \|f-v\|_{_\cT} 
  =\min_{{\scriptsize\btheta}\in\R^{M(n,d)}}
  \|f(\cdot)-v(\cdot;{\small\btheta})\|_{_\cT},
  \end{equation}
  
\begin{theorem}\label{error-bound1}
Assume that there exists a positive constant $\alpha$ such that $\alpha\, \|v\|^2 \leq \|v\|_{_\cT}^2$ for all $v\in \cM^1_{2n}(\sigma,d)$.
Let $f_{_\cT}$ be a solution of {\em (\ref{L2App-d})}.
Then there exists a positive constant $C$ such that
\begin{equation}\label{error-bound}
 \quad  \qquad C\,\|f-f_{_\cT}\|
 \leq \!\! \inf_{v\in \cM^1_{2n}(\sigma,d)}\!\! \left\{\|f-v\| + \!\!\sup _{w\in \cM^1_{2n}(\sigma,d)}\!\! \dfrac{|(\cI-\cQ_{_\cT})(vw)|}{\|w\|}\right\}+ \!\!\sup _{w\in \cM^1_{2n}(\sigma,d)}\!\! \dfrac{|(\cI-\cQ_{_\cT})(fw)|}{\|w\|}.
\end{equation}
\end{theorem}

\begin{proof}
Since ${\cal M}_n(\sigma,d)$ is a set, $f_{_\cT}\in {\cal M}_n(\sigma,d)$ is then characterized by the inequality
\begin{equation}\label{characterization}
    (f-f_{_\cT}, v-f_{_\cT})_{_\cT} \leq 0
    \quad\forall\,\, v\in {\cal M}_n(\sigma,d).
\end{equation}
For any $v\in {\cal M}_n(\sigma,d)$, it follows from the assumption and (\ref{characterization}) that
\begin{eqnarray*}
\alpha \|f_{_\cT}-v\|^2
&\leq & \|f_{_\cT}-v\|^2_{_\cT}
 %= (f_{_\cT}, f_{_\cT}-v)_{_\cT} -(v, f_{_\cT}-v)_{_\cT}
 \leq (f, f_{_\cT}-v)_{_\cT} -(v, f_{_\cT}-v)_{_\cT} \\[2mm]
 &=& \Big((f, f_{_\cT}-v)_{_\cT} -(f, f_{_\cT}-v)\Big) + \Big((v, f_{_\cT}-v) -(v, f_{_\cT}-v)_{_\cT}\Big) + (f-v,f_{_\cT}-v)
 \end{eqnarray*}
which, together with the Cauchy-Schwarz inequality, implies
 \[
 \alpha \|f_{_\cT}-v\|
 \leq  \!\!\sup _{w\in \cM^1_{2n}(\sigma,d)}\!\! \dfrac{|(\cI-\cQ_{_\cT})(fw)|}{\|w\|} + \!\!\sup _{w\in \cM^1_{2n}(\sigma,d)}\!\! \dfrac{|(\cI-\cQ_{_\cT})(vw)|}{\|w\|}+ \|f-v\|.
 \]
Combining the above inequality with the triangle inequality
 \[
 \|f-f_{_\cT}\| \leq \|f-v\| + \|v-f_{_\cT}\|
 \]
and taking the infimum over all $v\in \cM^1_{2n}(\sigma,d)$ yield (\ref{error-bound}). This completes the proof of the theorem.
\end{proof}

Theorem~\ref{error-bound1} indicates that the total error of the best least-squares approximation with numerical integration is bounded by the approximation error of the neural network and the error of the numerical integration. To ensure the approximation accuracy of the given neural network, we need to choose a numerical integration with a compatible accuracy, e.g., the composite ``mid-point'' numerical integration on an adaptively refined uniform partition. %Moreover, since $f^*,\,\, \hat{f} \in {\cal M}_n(\sigma,d)$ are piecewise linear functions, with a proper numerical integration, the second term on the left-hand side of (\ref{error-bound}) either vanishes or is of high order.

\begin{remark}
The assumption in {\em Theorem~4.1} is known as the Marcinkiewicz problem in literature and has not been verified for functions in $\cM^1_{2n}(\sigma,d)$. Recently, Temlyakov {\em \cite{Temlyakov2018}} introduced a new technique to systematically study this and related issues for functions in various finite dimensional subspaces. % but not yet for functions in nonlinear manifolds such as $\cM^1_{2n}(\sigma,d)$.
\end{remark}

\section{Adaptive network Enhancement (ANE) Method}

For a given target function $f(\bx)$, let $f_{_\cT}(\bx,{\small\btheta}^*_{_\cT})$ be the solution of problem (\ref{L2App-d}). For a given tolerance $\epsilon>0$, this section studies self-adaptive method for creating a two-layer ReLU NN and a numerical integration mesh such that the approximation accuracy is within the prescribed tolerance, i.e., 
 \begin{equation}\label{tolerance}
     \|f-f_{_\cT}\| \leq \epsilon \, \|f\|.
 \end{equation}

First, we consider the case that the numerical integration based on a partition $\cT$ is sufficiently accurate. Similar to the idea of the standard adaptive mesh-based numerical methods, we start with a two-layer ReLU NN with a small number of neurons, solve the optimization problem in (\ref{L2App-d}), and estimate the total error
by computing {\it a posteriori} error estimator
 \begin{equation}\label{estimator-NN}
     \xi = \|f-f_{_\cT}\|_{_\cT}/ \|f\|.
 \end{equation}
If $\xi >\epsilon$, we then enhance the NN by adding new neurons and this procedure repeats until (\ref{tolerance}) is met. This process is referred as the adaptive network enhancement (ANE) and it generates a two-layer ReLU NN whose approximation to $f$ satisfies a given approximation accuracy target. 

An immediate key question for the ANE method is: how many new neurons will be added at each adaptive step? %and how to initialize them? 
To address this issue, we propose two network enhancement strategies. One is global and the other is local. The global one is based on the assumption that the network approximation to the target function $f$ has a fixed convergence rate $\alpha$:
\[
\hat{\xi}^{(k)}= \|f-f^{(k)}_{_\cT}\|_{_\cT} = {\cal O}( n_{k}^{-\alpha}),
\]
where $f^{(k)}_{_\cT}$ is the approximation in ${\cal M}_{n_k}(\sigma,d)$, $n_k$ is the number of neurons of the $k^{th}$ NN, and $\alpha$ is the order of  approximation. A simple calculation suggests the following number of neurons for the next network:
\begin{equation}\label{global-n}
n_k=\min\left\{2n_{k-1}, \left\lceil \left(\hat{\xi}^{(k-1)}/\epsilon\right)^{1/\alpha_k}n_{k-1} \right\rceil\right\},
%\left\{\begin{array}{ll}\min\left\{2n_1, \lceil \hat{\xi}^{(1)}n_1/\epsilon \rceil\right\}, & k=2,\\[4mm]\min\left\{2n_{k-1}, \left\lceil \left(\hat{\xi}^{(k-1)}/\epsilon\right)^{\alpha_k}n_{k-1} \right\rceil\right\}, & k\ge 3\end{array}\right.
\end{equation}
where $\alpha_k$ is an approximation to the order $\alpha$. For $k\ge 3$,  $\alpha_k=\ln \left(\hat{\xi}^{(k-2)}\big/\hat{\xi}^{(k-1)}\right)\Big/\ln \left(n_{k-1}/n_{k-2}\right)$.
Possible choice for $\alpha_2$ is $1$ (linear rate) or some positive real number based on some {\it a priori} information of the target function.

To introduce our local network enhancement strategy, we notice that ${\cal M}_n(\sigma,d)$ is the set of continuous piece-wise linear functions with $n$ free hyper-planes given by (\ref{planes}). 
For any bounded domain $\Omega\in\R^d$, these $n$ hyper-planes plus the boundary of the domain $\Omega$ form a partition, $\cK_n=\{K\}$, of the domain $\Omega$. Again, the partition means that union of all subdomains of ${\cal K}_n$ equals the whole domain $\Omega$ and that any two distinct subdomains of ${\cal K}_n$ have no intersection.
We will refer to $\cK_n=\{K\}$ as the physical partition of the domain $\Omega$.

%These breaking hyper-planes become breaking points $x=b_i$ in one dimension and breaking lines $(\cos \bomega_i)\, x_1 + (\sin \bomega_i)\, x_2 = b_i$ in two dimensions. Therefore, each element $K\in \cK_n$ is a simple sub-interval in one dimension and a complicated polygonal sub-domain in $d\ge 2$ dimensions. 

This observation implies that the network enhancement strategy could make use of local errors on elements of the physical partition $\cK_n$. Specifically, let us introduce local error indicator $\xi_{_K}$ for each element $K\in\cK_n$:
 \begin{equation}\label{neuron-indicator}
    \xi_{_K}= \|f-f_{_\cT}\|_{_{K,\cT}}
    \equiv \left(\sum_{x_{K^\prime}\in K}(f-f_{_\cT})^2(\bx_{_{K^\prime}}) |K^\prime|
    \right)^{1/2}.
\end{equation}
%where $\hat{F}(\bx,\,\btheta)$ is the solution of problem (\ref{L2App-d}). % associated with the integration partition $\cT$. 
We then define a subset $\hat{\cK}_n$ of $\cK_n$ by using either the following average marking strategy:
 \begin{equation}\label{BV-marking-2}
    \hat{\cK}_n =\left\{K\in \cK_n\, :\,
    \xi_{_K}
    \ge \, \dfrac{1}{\#\cK_n}\sum_{K\in \cK_n}\xi_{_K}\right\},
\end{equation}
%with $\gamma_1\in (0,\,1)$
where $\#\cK_n$ is the number of elements of $\cK_n$, or the bulk marking strategy: finding a minimal subset $\hat{\cK}_n$ of $\cK_n$ such that
\begin{equation}\label{marking-2}
    \sum_{K\in \hat{\cK}_n} \xi^2_{_K}
    \ge \gamma_1\, \sum_{K\in \cK_n} \xi^2_{_K}
    \quad\mbox{for }\,\, \gamma_1\in (0,\,1).
\end{equation}
With the subset $\hat{\cK}_n$, the number of new neurons to be added to the NN is equal to the number of elements in $\hat{\cK}_n$.

With an accurate numerical integration, the ANE method is defined in Algorithm 5.1.
\begin{algorithm}
{\bf {\sc \bf Algorithm 5.1}} Adaptive two-layer ReLU NN with a fixed $\cT$.\\
Given a target function $f(\bx)$ and a tolerance $\epsilon >0$, starting with a two-layer ReLU NN with a small number of neurons, %For $k=0,\,1,\, ...$
\begin{itemize}
    \item[(1)] solve the optimization problem in (\ref{L2App-d});
    \item[(2)] estimate the total error by computing $\xi= \left(\sum\limits_{K\in \cK} \xi^2_{_K}\right)^{1/2}\, / \|f\|_{_\cT}$, where $\cK$ is the physical partition of the current approximation;
    \item[(3)] if $\xi< \epsilon$, then stop; otherwise, go to Step (4);
    \item[(4)] add new neurons to the network by using the network enhancement strategy, then go to Step (1). 
\end{itemize}
\end{algorithm}

\smallskip

Next, we consider adaptive mesh refinement (AMR) on numerical integration for a fixed NN.
Let $f_{_\cT}(\bx,\,{\small\btheta}^*_{_\cT})$ be the solution of problem (\ref{L2App-d}) associated with the partition $\cT$.  
Let $\hat{\cT}$ be a subset of $\cT$ generated by using either the average or the bulk marking strategy.
For each marked element $T\in\hat{\cT}$, this $d$-dimensional cube is subdivided into $2^d$ small cubes of equal size. The new partition $\cT^\prime$ consists of elements in $\cT\setminus \hat{\cT}$ and new elements generated from $\hat{\cT}$. Denote by $f_{_{\cT^\prime}}(\bx,{\small\btheta}^*_{_{\cT^\prime}})$ the solution of problem (\ref{L2App-d}) associated with the partition $\cT^\prime$. For both solutions $f_{_{\cT}}$ and $f_{_{\cT^\prime}}$
based on the mesh $\cT$ and its refinement $\cT^\prime$, define the following global estimators:
 \[
 \eta (f_{_\cT})= \left(\sum\limits_{T\in \cT}\eta_{_{T}}(f_{_\cT})^2 \right)^{1/2}
 \quad\mbox{and}\quad
 \eta (f_{_{\cT^\prime}})= \left(\sum\limits_{T\in \cT^\prime}\eta_{_{T}}(f_{_{\cT^\prime}})\right)^{1/2}.
 \]
where local indicators on $\cT^\prime$ are given by
\begin{equation}\label{Q-indicator}
\eta_{_{T}}(f_{_\cT})= \|f-f_{_{\cT}}\|_{_{T,\cT}}
    %\equiv \big|(f-\hat{F}_{_\cT})(\bx_{_T})\big|
    \quad\mbox{and}\quad
    \eta_{_{T}}(f_{_{\cT^\prime}})= \|f-f_{_{\cT^\prime}}\|_{_{T,\cT^\prime}}.
    %\equiv \big|(f-\hat{F}_{_\cT})(\bx_{_T})\big|
    %\quad\mbox{for all }\,\, T\in\cT^\prime,
\end{equation}

The following algorithm generates a numerical integration mesh which ensures approximation accuracy of a given NN.
\begin{algorithm}
{\bf {\sc \bf Algorithm 5.2}} Adaptive Mesh Refinement with a fixed NN.\\
Given a target function $f(\bx)$ and  
the solution of problem (\ref{L2App-d})
on the partition $\cT$,
\begin{itemize}
    %\item[(1)] compute local error indicator 
    %$\eta_{_T}$ for all $T\in\cT$; %by computing $\xi= \left(\sum\limits_{K\in \cK} \xi^2_{_K}\right)^{1/2}$;
    \item[(1)] refine $\cT$ by the refinement strategy to obtain a new partition $\cT^\prime$
    \item[(2)] solve the minimization problem in (\ref{L2App-d}) on $\cT^\prime$;
    \item[(3)] if $\eta (f_{_{\cT^\prime}}) \leq \gamma_{_2} \eta (f_{_{\cT}})$, then %replace $\cT^\prime$ by $\cT$ and 
    go to Step (1) with $\cT = \cT^\prime$; otherwise, output $\cT$.
\end{itemize}
\end{algorithm}

\smallskip

The stopping criterion used in Algorithm~5.2 is based on whether or not the mesh refinement on numerical integration improves approximation accuracy.  When the refinement does not improve accuracy much, the AMR stops and outputs the current mesh.

Finally, we are ready to present adaptive network enhancement (ANE) method for a two-layer ReLU NN including AMR for numerical integration in Algorithm 5.3. The purpose of the AMR for numerical integration is to ensure approximation accuracy with less quadrature points than a fine uniform partition. Comparing with the ANE, the AMR is secondary. 

\smallskip

\begin{algorithm}
{\bf {\sc \bf Algorithm 5.3}} Adaptive two-layer ReLU NN.\\
Given a target function $f(\bx)$ and a tolerance $\epsilon >0$, starting with a coarse uniform partition $\cT_0$ of the domain $\Omega$ for numerical integration and with a two-layer ReLU NN with a small number of neurons, 
\begin{itemize}
    \item[(1)] solve the minimization problem in (\ref{L2App-d}); %with $n_k$ neurons on the numerical partition $\cT_k$;
    \item[(2)] use Algorithm~5.2 to generate a numerical integration mesh $\cT$;
    \item[(3)] solve the minimization problem in (\ref{L2App-d}) associated with $\cT$;
    \item[(4)] estimate the total error by computing $\xi= \left(\sum\limits_{K\in \cK} \xi^2_{_K}\right)^{1/2} /\,\|f\|_{_\cT}$, where $\cK$ is the physical partition of the current approximation;
    \item[(5)] if $\xi< \epsilon$, then stop; otherwise, go to Step (6);
    \item[(6)] add new neurons to the network by using the network enhancement strategy, then go to Step (1).
\end{itemize}
\end{algorithm}

\smallskip

\section{Strategies for training (iterative solvers)}

%The minimization problem in (\ref{L2App}) is often 
%solved numerically by iterative methods such as gradient descent (GD), Stochastic GD, Adam, etc. (see, e.g., \cite{BoCuNo2018} for a review paper in 2018). Since the loss function $\|F(\btheta)-f\|^2_{_\cT}$ is a non-convex function of parameters $\btheta$, initial of any iterative method is crucial for finding a ``good'' approximation. 
The exceptional power of DNNs in approximation come with a price: the procedure for determining the values of the parameters is now a problem in nonlinear optimization. This high dimensional, nonlinear optimization problem tends to be computationally intensive and complicated. Currently, it is often solved by iterative optimization methods such as gradient descent (GD), Stochastic GD, Adam, etc. (see, e.g., \cite{BoCuNo2018} for a review paper in 2018 and references therein). 
Usually nonlinear optimizations have many solutions, and the desired one is obtained only if we start from a close enough first approximation. The ANE method provides a natural process for obtaining a good initialization. This section describes our initialization for all three stages of the ANE method.

%Starting with a relatively small NN, the approximation of the previous NN is already a good approximation to the current NN in the loops of the ANE method. To provide a better approximation than the previous one, we divide all network parameters into two groups: linear parameters (output layer weights and bias) and nonlinear parameters (hidden layer weights and biases). Initialization of nonlinear parameters are based on their physical partitioning of the domain and initial of linear parameters are obtained by solving a system of linear equations with given nonlinear parameters. 

The first stage is the beginning of the ANE method, in which we specify the size of the NN, both input and output weights and bias, and a partition of the domain for numerical integration. Due to the fact that input weights and bias determine physical locations of breaking hyper-planes, we first subdivide the domain $\Omega$ by a coarse, uniform partition and then distribute those breaking hyper-planes on the mesh of this partition. For example, when $\Omega=(0,\,1)^2$, the two-layer NN with $2(m_0+1)$ neurons use the following initial breaking lines: 
 \[
 x= \dfrac{i}{m_0}
 \quad\mbox{and}\quad
 y= \dfrac{i}{m_0}
 \quad\mbox{for }\,\, i=0,\,1,\, ..., \, m_0.
 \]
This breaking lines imply the following input weights and bias:
\[
{\small\btheta}_{1,i} = \left( (1,\,0),\, \dfrac{i}{m_0}\right)
\quad\mbox{and}\quad
{\small\btheta}_{2,i} = \left( (0,\,1),\, \dfrac{i}{m_0}\right)
\]
for $i=0,\,1,\, ..., \, m_0$.
For numerical integration, we again start with a uniform partition $\cT$ of the domain $\Omega$
which, in general, is much finer than the previous physical partition initializing the NN. Initial of the output weights and bias is given by the solution of the system of linear equations in (\ref{O_W}).

The second stage is the AMR for numerical integration. For each new partition $\cT$, natural initial of parameters ${\small\btheta}$ is the corresponding values of the current approximation since the NN remains unchanged. 

The third stage is when the NN is enhanced by adding new neurons. Clearly, parameters corresponding to old neurons will use the current approximation as their initial. To initialize corresponding parameters of new neurons, for the global enhancement strategy, one can add new neurons randomly; or add new neurons uniformly across the domain (i.e. set their input weights and biases with corresponding break hyper planes uniformly subdividing the domain).  For the local enhancement strategy, we propose to make use of the subset $\hat{\cK}_n$ marked in (\ref{BV-marking-2}) or (\ref{marking-2}). For each element $K\in\hat{\cK}_n$, we add one neuron whose initial is corresponding to the breaking hyper-plane that passes through the centroid of $K$ and orthogonal to the direction vector with the smallest variance of quadrature points in $K$. This direction vector may be computed by the Principal Component Analysis method (or PCA \cite{PCA1901}).  %There are two ways to initialize 
For output weights and biases corresponding to new neurons, a simple initial is to set them zero. This means that the initial of the approximation is the current approximation. A better way is to solve problem (\ref{O_W}) for all output weights and bias by using the current breaking hyper-planes for the input weights and bias.

\section{Numerical Experiments}
In this section, we present our numerical experiments on using ANE to approximate various functions. In all experiments, the minimization problem  (\ref{L2App-d}) is solved using the Adam version of gradient descent \cite{kingma2014adam}.  For each run during the adaptive process, the stopping criteria for the iterative solver is set as follows: the solver stops when the loss function $\|f-\hat{f}\|_{_\cT}$ decreases within $0.1\%$ in the last 2000 iterations. This stopping criteria is set to explore the network approximation power without constraining the number of iterations.

\subsection{Smooth Function}
The first test problem is a smooth function of one variable
\begin{equation}\label{test1}
 f(x)=x\left(e^{-(x-\frac{1}{3})^2/k}-e^{-\frac{4}{9}/k}\right), 
\end{equation}
which is defined on the interval $\Omega = [0,1]$. When $k=0.01$, this function is the solution to a Poisson equation studied in  \cite{he2018relu,cai2020}. We use this simple toy problem to test the efficacy of the proposed ANE method.

\begin{figure}[htbp]
  \centering 
   \subfigure[Initial NN model with 10 uniform break points]{ 
    \label{figpoisson:a} 
    \includegraphics[width=2.4in]{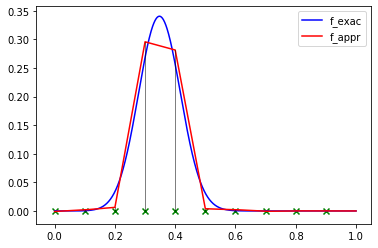}} 
  \hspace{0.3in} 
  \subfigure[Optimized NN model with 10 neurons]{ 
    \label{figpoisson:b} 
    \includegraphics[width=2.4in]{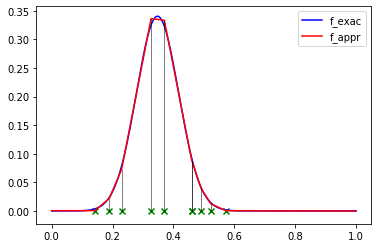}} 
     \hspace{0.3in} 
    \subfigure[Optimized NN model with 13 neurons using ANE]{ 
    \label{figpoisson:c} 
    \includegraphics[width= 2.4in]{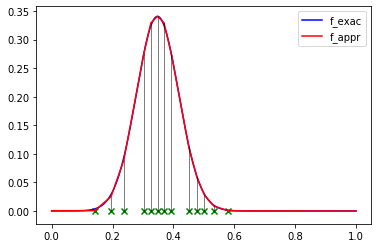}} 
     \hspace{0.3in} 
    \subfigure[Optimized NN model with 20 neurons using ANE]{ 
    \label{figpoisson:d} 
    \includegraphics[width= 2.4in]{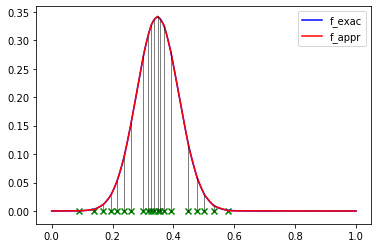}} 
     \hspace{0.3in} 
    \subfigure[Fixed NN model with 20 neurons]{ 
    \label{figpoisson:e} 
    \includegraphics[width= 2.4in]{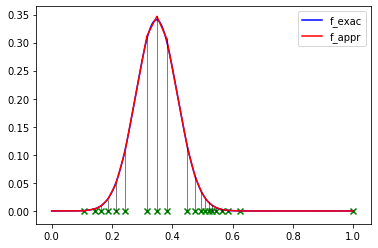}} 
     \hspace{0.3in} 
    \subfigure[Fixed NN model with 38 neurons]{ 
    \label{figpoisson:f} 
    \includegraphics[width= 2.4in]{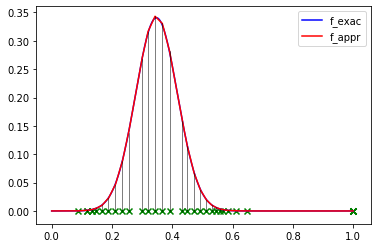}} 
     \hspace{0.3in}     
  \caption{Results of using two-layer ReLU networks for approximating function (\ref{test1})}
  \label{poisson_app} %% label for entire figure 
\end{figure}

The target approximation accuracy is set as $\epsilon=0.005$. A fixed uniform partition $\cT$ with $1000$ quadrature points is used for this experiment. We start from $10$ neurons for the input layer with their break points initialized uniformly across the domain, i.e., $b_i=0.1i$ for $i=0,1,\cdots,9$. The initial network model's output weights and biases are set by solving the linear system in (\ref{O_W}). This initial model is shown in Fig.\ref{figpoisson:a}. 

After the first run network training (solving (\ref{L2App-d})  using the Adam solver), the network adjusts its parameters to adapt the target function $f$. The resulting optimized network model with $10$ neurons is shown in Fig.\ref{figpoisson:b}. This NN model provides a near-optimum free-knot piecewise linear spline with $10$ knots shown as the green break points in the Figure. Base on the partition of the domain with current set of break points, we adopt the average marking strategy (\ref{BV-marking-2}) to mark the elements with errors larger than average, and then add neurons accordingly by setting the newly added neuron's initial biases at the centers of the elements to be refined. We then resolve for the new output layer's parameters using (\ref{O_W}) and trained the network for the second run. This process repeats until the approximation error is lower than the target $\epsilon$. The ANE method iterates itself three runs from $10$ to $13$ then to $20$ neurons. The intermediate result at $13$ neurons is depicted in Fig.\ref{figpoisson:c}. The ANE process ends at $20$ neurons, which gives a relative approximation accuracy of $\xi=0.003837$, falling below the target $\epsilon$. 

\begin{figure}[htbp]
  \centering 
   \subfigure[10 neurons (with three marked elements)]{ 
    \label{Poisson_err_dist:a} 
    \includegraphics[width=1.8in]{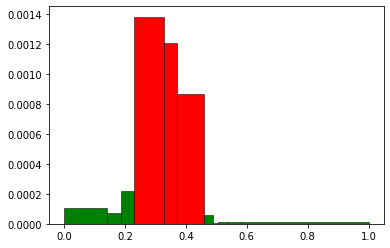}} \hfill 
   \subfigure[13 neurons (with seven marked elements]{ 
    \label{Poisson_err_dist:b} 
    \includegraphics[width=1.8in]{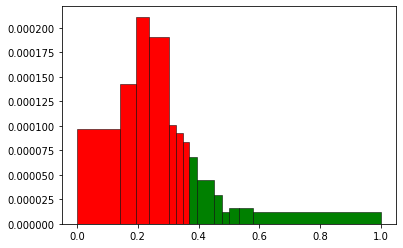}} \hfill   
    \subfigure[20 neurons]{ 
    \label{Poisson_err_dist:c} 
    \includegraphics[width=1.8in]{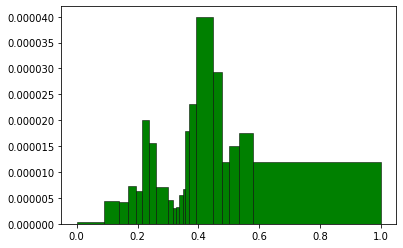}} 
  \caption{Error distribution on physical partitions generated in the ANE process for the first test problem, where red partitions are the elements to be refined. }
  \label{Poisson_err_dist} %% label for entire figure 
\end{figure}

In this one-dimensional problem, we utilize a fixed learning rate of 0.001. Fig. \ref{Poisson_err_dist:a} and \ref{Poisson_err_dist:b} show the error per element distribution on the physical partition generated through the iterative process. In those two figures, red bars correspond with marked elements where new neurons are to be added.  Marked red elements are refined and the iterative process gradually drags all elements error down to a smaller scale with a trend to distribute the error evenly among the physical partitions, see error distribution of the final network model in Fig. \ref{Poisson_err_dist:c}.

%while the green bars corresponds lower error elements. Elements with red bar errors are refined and iterative process gradually drags all elements error down to a smaller scale while tending to distribute them evenly.
%This is based on the following observation: errors distributed on the physical partition (mesh) from a trained network model can be divided into two groups, elements with lowers errors and larger errors, and there is gap between them. Lower error elements correspond to nearly linear regions of $f$, where piece-wise linear function can approximate them accurately. While the larger error elements are associated with non-linear zones of $f$. One guess for larger error elements is that, the optimizer tends to move the break points such that the errors distributed to the non-linear zones on each element are nearly uniform. This trends is verified in our experiment. 

\begin{figure}[htbp]
\centering
\includegraphics[width=3in]{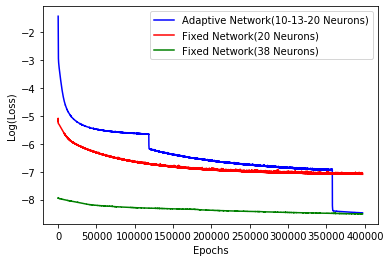}
\caption{Log training loss with three different network models in the first numerical experiment.}
\label{poisson_loss}
\end{figure}

We further compare the performance of our adaptive network structure with a network model of fixed number of neurons. This is to check if the adaptive process has a potential to land in a better global minimum. The comparison results are illustrated in Table \ref{poissonnumerical2} and Fig. \ref{poisson_loss}. In particular, they show that approximation accuracy of the adaptive network using ANE with $20$ neurons is almost same as that of the fixed network with $38$ neurons, and is better than the fixed network structure of the same size. This experiment indicates that fix networks might tend to be trapped in local minimums. 

Finally, we test the performance of enhancement strategy using (\ref{global-n}), and compare two methods of initialization under the global adaptive enhancement scheme. The first initialization method is to add new neurons randomly and set their corresponding output weights as zeros as initial; and the second method is to add new neurons uniformly across the domain and solve (3.4) for output weights and bias.  Table \ref{poissonnumerical3} list the results of this experiment. Due to the non-convex optimization, one can see that different initialization strategy result in differences in approximation results. The first initialization method is easily trapped in local minimum under the Adam optimizer. Considering the randomness in the initial of newly added neurons, we repeat this test three times and report the best result in the table. While in the second strategy, we start from a better point using global uniform refinement, this results in a better performance. However, the uniform initial strategy during the adaptive process does not consider the error distribution evaluated from the previous stage, which explains why it is still inferior to the local error based marking and refinement strategy. For the rest experiments, we only use the local enhancement method.

\begin{table}[htb]
\caption{Comparing adaptive neural network with fixed networks for testing problem {\em (\ref{test1})}}
\label{poissonnumerical2}
\centering
\begin{tabular}{  |l |c |c | c | c |  c| c }
	\hline
	Network (neurons) & \# Parameters & $\|f-f_{_\cT}\|_{_\cT}/ \|f\|$  \\[2mm] \hline
	Fixed (20)  & 41 & 0.007644     \\ \hline 
	Fixed (38)  & 77 & 0.003762   \\ \hline 
	Adaptive (10$\rightarrow$13$\rightarrow$20) & 41 & 0.003837  \\ \hline
\end{tabular}
\end{table}

\begin{table}[htb]
\caption{Global network enhancement and initialization strategy for tesing problem {\em (\ref{test1})}}
\label{poissonnumerical3}
\centering
\begin{tabular}{  |l |c |c |c | c | c |  c| c }
	\hline
	Network (adaptive neurons) & \ Initialization &\# Parameters & $\|f-f_{_\cT}\|_{_\cT}/ \|f\|$  \\[2mm] \hline
	10$\rightarrow$20$\rightarrow$40  & random & 81 & 0.005221    \\ \hline 
	 10$\rightarrow$20$\rightarrow$30 & uniform & 61 & 0.004455  \\ \hline
\end{tabular}
\end{table}

\subsection{Functions with intersecting interface singularities}
\begin{figure}[htb!]
  \centering 
   \subfigure[Target function f(x,y)]{ 
    \label{kellog:a} 
    \includegraphics[width=1.9in]{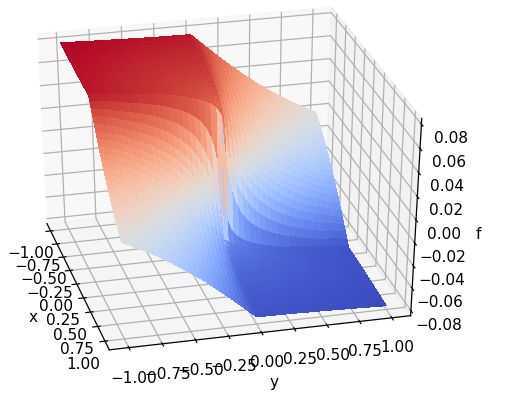}} \hfill 
    \subfigure[Initial break lines \newline(20 neurons, 100 elements)]{ 
    \label{kellog:b} 
    \includegraphics[width=1.7in]{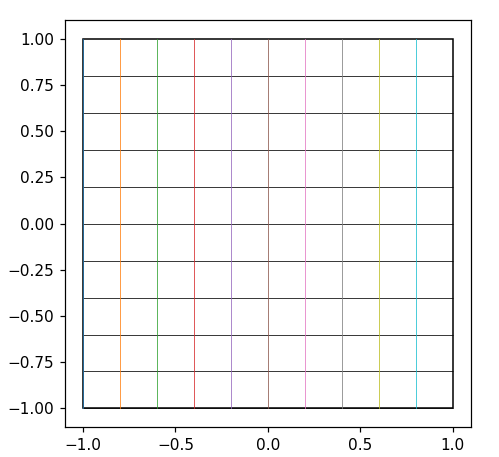}} \hfill 
   \subfigure[Initial NN model with 20 neurons]{ 
    \label{kellog:c} 
    \includegraphics[width=1.9in]{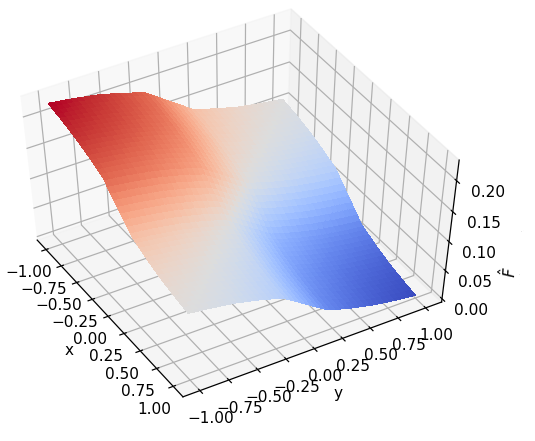}} \hfill   
    \vspace{0.2in} 
    \subfigure[Optimum break lines \newline (20 neurons, 103 elements)]{ 
    \label{kellog:d} 
    \includegraphics[width=1.8in]{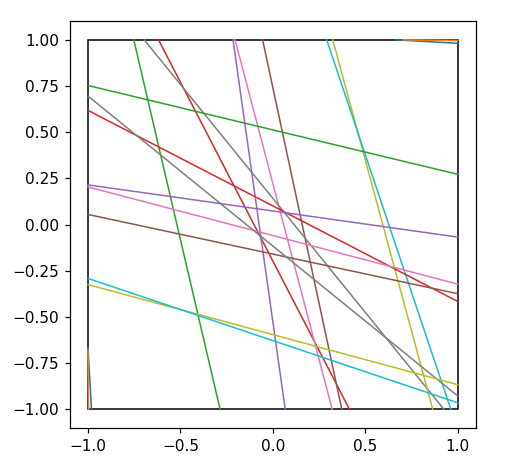}} \hfill     
    \subfigure[Optimum NN model of 20 neurons,  $\xi=0.038733$]{
    \label{kellog:e} 
    \includegraphics[width=1.9in]{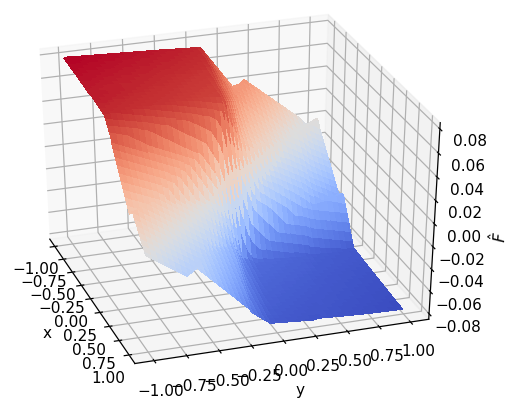}}\hfill 
    \subfigure[15 more neurons added with break lines in light blue]{ 
    \label{kellog:f} 
    \includegraphics[width=1.7in]{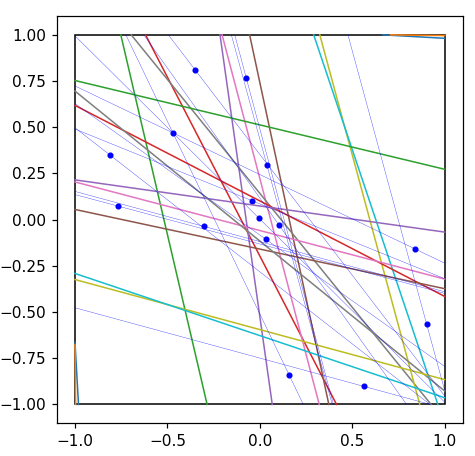}} \hfill
    \vspace{0.2in} 
    \subfigure[Optimum break lines of 35  \newline neurons and 34 more neurons are to be added in the second run]{ 
    \label{kellog:g} 
    \includegraphics[width=1.7in]{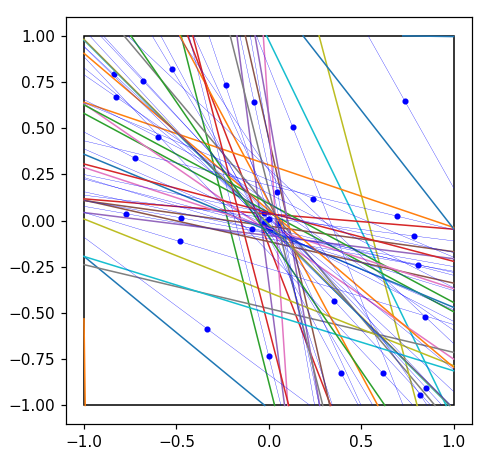}}\hfill  
    \subfigure[Optimum break lines \newline (69 neurons, 1286 elements)]{ 
    \label{kellog:h} 
    \includegraphics[width=1.7in]{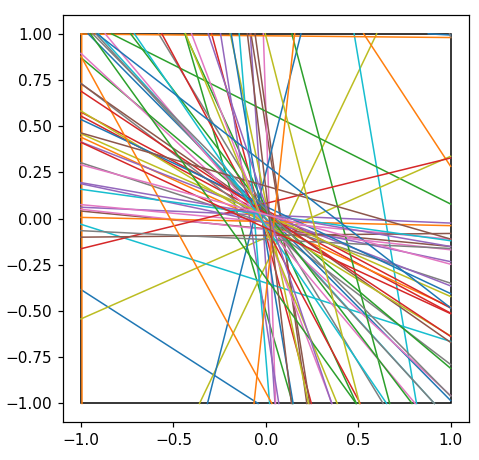}}\hfill  
    \subfigure[Optimum NN model of 69 \newline neurons, $\xi=0.008476$]{ 
    \label{kellog:i} 
    \includegraphics[width=2in]{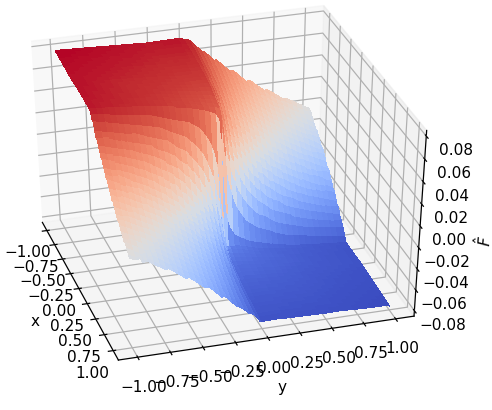}}\hfill  
  \caption{ANE results of using 2-layer ReLU networks for approximating function in {\em (\ref{test2-1})}}
  \label{singular} %% label for entire figure 
\end{figure}

This section reports the numerical results for a two-dimension problem with intersecting interface singularity. Let $\Omega = (-1,1)^2$ and 
\begin{equation}\label{test2-1}
f(r, \theta) = r^{\beta} \mu(\theta)
 \end{equation}
in the polar coordinates at the origin with

%\begin{equation}\label{test2-2}
$$
 \mu(\theta)=
\left\{ \begin{array}{rclll}
\cos((\pi/2-\sigma)\beta)\cdot \cos((\theta-\pi/2+ \rho)\beta), \quad &\text{if} & 0\leq \theta\leq\pi/2,\\[1mm]
\cos(\rho\beta)\cdot \cos((\theta-\pi+ \sigma)\beta), \quad &\text{if} & \pi/2 \leq \theta\leq\pi,\\[1mm]
\cos(\theta\beta)\cdot \cos((\theta-\pi-\rho)\beta), \quad &\text{if} & \pi \leq \theta\leq 3\pi/2,\\[1mm]
\cos((\pi/2-\rho)\beta)\cdot \cos((\theta-3\pi/2-\sigma)\beta), \quad &\text{if} & 3\pi/2 \leq \theta\leq 2\pi,
 \end{array}\right.
$$ 
 %\end{equation}
where $\beta =0.1$, $\sigma=−14.92256510455152$, and $\rho=\pi/4$ are parameters. The function $f(r, \theta)$ (see Fig.\ref{kellog:a}) is the solution of the elliptic interface problem with intersecting interface singularity and a benchmark test problem for adaptive finite element method (see,  e.g., \cite{Morin2002, Cai2009}).
 
We test the ANE method with a fixed integration mesh using $400\times 400$ quadrature points. %while the network's approximation performance was evaluated using a finer mesh of $1000\times 1000$ quadrature points. 
The target approximation accuracy is set as $\epsilon=0.01$.  The ANE process starts with a small network of $20$ neurons, and the network is initialized such that the break lines are distributed evenly in the domain, with half of them parallel to $x$-axis ($\bomega_i=0$ and $b_i=-1+0.2i$ for $i=0,\cdots,9$) and the other half parallel to $y$-axis ($\bomega_i=\pi/2$ and $b_i=-1+0.2(i-10)$ for $i=10,\cdots,19$). See Fig. \ref{kellog:b} for the initial partition of the domain. The initial network model using this uniform physical partition is obtained by solving the linear system in (\ref{O_W}) and is shown in  Fig. \ref{kellog:c}. After the first run network training, the optimum break lines corresponding to the $20$-neuron two-layer ReLU network is shown in Fig. \ref{kellog:d} and the corresponding network model is plotted in Fig. \ref{kellog:e}. With $20$ neurons ($61$ parameters), the adaptive network can approximate the target function $f$ in (\ref{test2-1}) with a relative error $\xi=0.038733$.
 
To achieve the target accuracy, the ANE calculates per element error base on the automatic generated physical partition of the domain.  Elements with relative large errors are marked using the bulk marking strategy (\ref{marking-2}) with $\gamma_1 =0.7$ (see the $15$ elements with blue dots shown in Fig \ref{kellog:f}). ANE process adds the same number of neurons as the $15$ marked elements, and those new neurons are initialized as follows: their corresponding breaking lines pass through the centroids of marked elements, with their directions aligned with the maximum principal directions of each geometric element. See Fig \ref{kellog:f} for the initial physical partition at the second run with the newly added neuron's breaking lines drawn in light blue. The second run network training converged at a relative error $\xi=0.019582$ (see the generated physical partition and marked elements in Fig \ref{kellog:g}). The ANE process stops at $69$ neurons with the corresponding physical partition and network model plotted in Fig \ref{kellog:h} and Fig \ref{kellog:i}. Notice that to calculate per element error, and to find an element's centroid and principal direction, we group the quadrature points located in the same element and use the point set within the element to compute its local error, centroid and PCA. This approximation method has an advantage of its computational simplicity; by avoiding calculation of the exact geometric shape of each element, this method can be easily extended to higher dimension problems or higher order activation functions. 

A fixed learning rate of $10^{-3}$ is adopted in this ANE process. The final network model achieves a $L^2$ relative error of $\xi=0.008476$, which meets our approximation accuracy target. The generated physical partition is highly adapted to the target function. Notice there is a point singularity around the origin in the function $f$, while the physical partition obtained in the adaptive network adjusts its elements shape and size such that the partition is dense around the singular point, this is a very favorable property of using NN model to approximate functions with singularities. Comparing with adaptive finite element methods (AFEMs) (see, e.g., \cite{Cai2009}), the ANE method has much fewer degrees of freedom than AFEMs. 

To evaluate the effect of numerical integration to the total approximation error, we tested a two-layer network of $69$ neurons using varying ${\cal T}$ with different number of quadrature points. 
The results are given in Table \ref{kellog_quad}. As shown in the table, with finer integration meshes of more number of quadrature points, the integration accuracy can be improved (refer to the `Integration accuracy' column in the table). Meanwhile, training a network model on finer mesh is harder which results in a lower training accuracy (see the `Training accuracy' column). However, the approximating power to the 
true function $f$ is improved (see the `Testing accuracy' column in Table \ref{kellog_quad}). Notice here the testing accuracy is estimated using a fine mesh ${\cal T^\prime}$ of $1000 \times 1000$ quadrature points. The gap between training accuracy and testing accuracy is reduced when more number of quadrature points is adopted.  This experiment also shows that the adaptive network may achieve better approximation result compared with the fixed network of the same size, see the last two rows in Table \ref{kellog_quad}.

\begin{table}[htb]
\caption{The effect of numerical integration for the second testing problem {\em (\ref{test2-1})}}
\label{kellog_quad}
\centering
\begin{tabular}{  |l |c |c | c | c | c |  c| c }
	\hline
  \makecell{Network\\($\#$ quadrature)} & \makecell{Integration accuracy\\ $|(\cI-\cQ)(f)|/|\cI(f)|$} & \makecell{Training accuracy\\ $\|f-f_{_\cT}\|_{_\cT}/ \|f\|$} & \makecell{Testing accuracy \\$\|f-f_{_\cT}\|_{_{\cT^\prime}}/ \|f\|$}  \\[4mm] \hline
    Fixed  (50x50)  &0.002638 & 0.007885 & 0.013187\\ \hline
	Fixed  (100x100)  &0.000753 & 0.008515 & 0.010257 \\ \hline 
	Fixed  (200x200)  &0.000462 & 0.009319  & 0.009877  \\ \hline 
	Fixed  (400x400)  &0.000370 & 0.009702  & 0.009850  \\ \hline
	ANE (400x400)   &0.000370 & 0.008319  & 0.008476  \\ \hline	
\end{tabular}
\end{table}

\subsection{Functions with transition layers}
\begin{figure}[htb!]
  \centering 
  \subfigure[The target function $f$ with a circular transition layer]{ 
    \label{figsingular:a} 
   \includegraphics[width=2.5in]{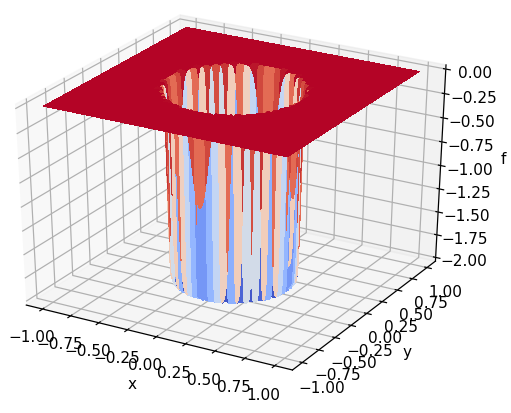}} 
    \hspace{0.2in}
  \subfigure[The generated AMR of 22201 quadrature points]{ 
    \label{figsingular:b} 
    \includegraphics[width=2.5in]{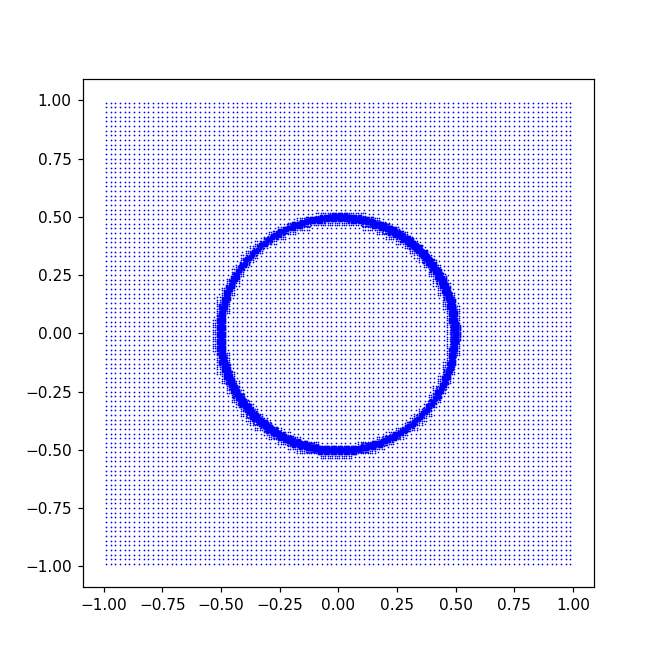}}
        \vspace{0.2in}
  \subfigure[ANE network model of 578 neurons using AMR in (b) for training (integration), $\xi=0.048771$]{ 
    \label{figsingular:c} 
    \includegraphics[width=2.5in]{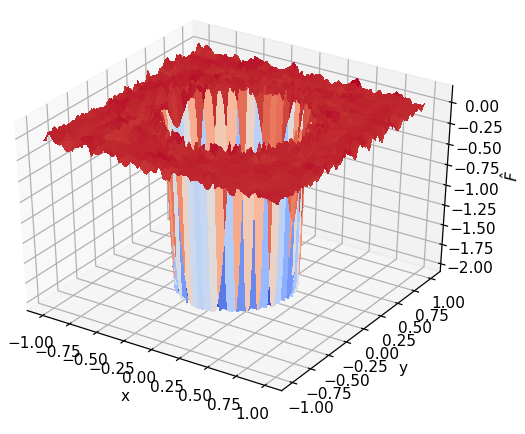}} 
        \hspace{0.2in}
  \subfigure[The physical partition generated by ANE of 578 neurons]{ 
    \label{figsingular:d} 
    \includegraphics[width=2.4in]{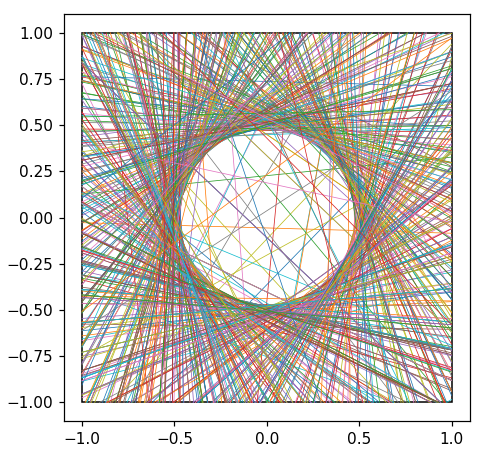}} 
      %\subfigure[a 3-layer network model to approximate target $f$ ]{ 
    %\label{figsingular:e} 
   %\includegraphics[width=2.5in]{Figures/singular2d_3layer(20neurons).png}} \hfill 
  \caption{ANE with AMR results of using 2-layer ReLU networks for approximating function in (\ref{test3}).}
  \label{singular-1} %% label for entire figure 
\end{figure}

The last problem we tested is a two-dimensional function with a transition layer around a circular region:
\begin{equation}\label{test3}
f(x,y) = \tanh\left(\frac{1}{\varepsilon}(x^2 +y^2 - \frac{1}{4})\right) - \tanh\left(\frac{3}{4\varepsilon}\right)
 \end{equation}
%\[f(x,y) = \tanh\left(\frac{1}{\varepsilon}(x^2 +y^2 - \frac{1}{4})\right) - %\tanh\left(\frac{3}{4\varepsilon}\right)
%\]
defined on the domain $\Omega = [-1,1]^2$. By varying $\varepsilon$, this type of functions show different level of difficulties due to the presence of transition layers. We set $\varepsilon=0.01$ in this experiment, and the corresponding function $f$ presents a large transition in a sharp circular zone, as shown in Fig. \ref{figsingular:a}. 

For this problem, we ran three tests to compare the results of using an uniform integration mesh versus adaptive mesh refinement (AMR). (1) The first experiment utilizes an uniformly allocated $400 \times 400 =1.6\times 10^5$ quadrature points and the ANE Algorithm 5.1 to obtain a network model of $578$ neurons with target accuracy $\epsilon=0.05$. (2) The second experiment uses Algorithm 5.3 which generates an AMR of $22201\approx 2.2\times 10^4$ quadrature points (as shown in Fig.\ref{figsingular:b}) and an adaptive NN of 578 neurons as well. The $22201$ quadrature points are generated by adaptive local mesh refinement of an initial mesh of $100\times 100$ quadrature points, using average marking strategy.  We set the last run ANE process to stop at $578$ neurons to allow a fair comparison to the first experiment.  (3) the third experiment matches the number of quadrature points used in the second experiment, but with those $150 \times 150$ quadrature points allocated uniformly across the domain, and a fixed network model of $578$ neurons was tested to compare the approximation performances with the ANE network using AMR integration mesh. 

The comparison results are illustrated in Table \ref{test3_result}. The ANE method using AMR for numerical integration achieves better performance compared with a finer uniform mesh of six times more quadrature points and it is superior compared with the similar mesh size but evenly distributed quadrature points. If limited computational resources are allocated which allows only certain number of quadrature points for numerical integration and network training, allocating quadrature points using AMR might achieve better approximation performance compared with the uniformly allocated quadrature points.

\begin{table}[htb!]
\caption{Networks approximation performances of uniform v.s. AMR integration mesh}
\label{test3_result}
\centering
\begin{tabular}{  |l |c |c | c | c | c | c |  c| c }
	\hline
	\makecell{Integration \\mesh}
  &\# quadrature &\# neurons & \makecell{Training accuracy\\$\|f-f_{_{\cT}}\|_{_\cT}/ \|f\|$} & \makecell{Testing accuracy \\$\|f-f_{_{\cT^\prime}}\|_{_{\cT^\prime}}/ \|f\|$}  \\[1mm] \hline
    Uniform &400x400 & ANE 578 & 0.050552 & 0.050587\\ \hline
	AMR &22201  &ANE 578 & 0.047423 & 0.048771 \\ \hline 
	Uniform  &150x150 &Fixed 578 &0.052497   & 0.053040  \\ \hline 
\end{tabular}
\end{table}

\begin{table}[htb!]
\caption{Approximation performances of a two-layer v.s. a three-layer NN}
\label{test3_result2}
\centering
\begin{tabular}{  |l |c |c | c | c | c | c |  c| c }
	\hline
	\makecell{NN structure\\(neurons)}
  &\#Quadrature &\#Parameters & \makecell{Training accuracy\\$\|f-f_{_{\cT}}\|_{_\cT}/ \|f\|$} & \makecell{Testing accuracy \\$\|f-f_{_{\cT^\prime}}\|_{_{\cT^\prime}}/ \|f\|$}  \\[1mm] \hline
	\makecell{Two-layer \\(578)} & AMR 22201  &1735 & 0.047423 & 0.048771 \\ \hline 
	\makecell{Three-layer\\ (20-20)}  &uniform 150x150 &501 &0.033751   & 0.033969\\ \hline 
\end{tabular}
\end{table}

 \begin{figure}[htbp]
  \centering 
    \subfigure[a three-layer network model of 20 neurons in each hidden layer, $\xi=0.033967$ ]{ 
    \label{figsingular2:a} 
   \includegraphics[width=2.8in]{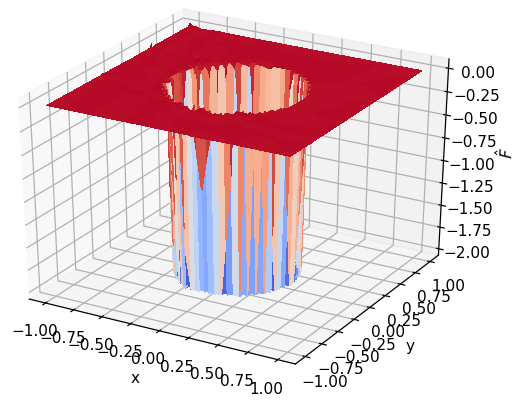}} \hfill 
       \hspace{0.2in}
    \subfigure[Physical partition generated by the three-layer network (black lines are the break lines in the first hidden layer, colored lines are the break polylines in the second hidden layer]{ 
    \label{figsingular2:b} 
   \includegraphics[width=2.6in]{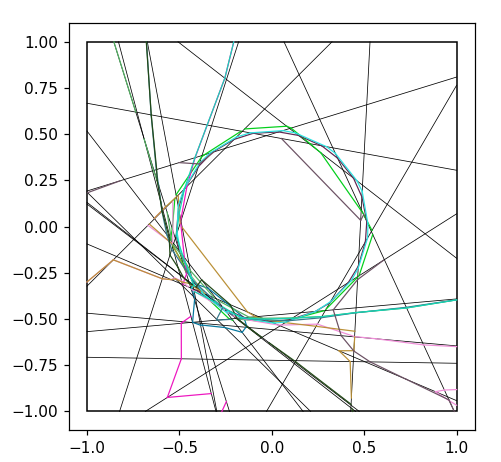}} \hfill 
  \caption{Approximation results of using a three-layer ReLU network for approximating function in (\ref{test3}).}
  \label{singular-2} %% label for entire figure 
\end{figure}

The function approximation result shown in Fig.\ref{figsingular:c} exhibits a certain level of oscillation which is not acceptable in some applications. Notice that the generated physical partition (see Fig.\ref{figsingular:d}) does capture the circular transition layers well when using $578$ break lines. However, this partition is too dense in the region where the function does not fluctuate much. A deeper ReLU network, which provides piece-wise breaking lines, might work better for this testing case. We verified this conjecture by using a three-layer ReLU network to approximate this function. Each hidden layer was set as fixed $20$ neurons which defines a network model of $501$ parameters. The relative approximation error $\xi$ using this three-layer ReLU network is $0.033967$. Comparing to the $578$ neurons and $1735$ parameters we used previously in the two-layer networks, a three-layer ReLU network of smaller size can approximate the same function with better accuracy (see Table. \ref{test3_result2}). As illustrated in Fig.\ref{figsingular2:a}, a three-layer network can reduce the oscillation exhibited in the shallow network, while archiving a better approximation accuracy with less complicated domain partition, see Fig.\ref{figsingular2:b} for the physical partition generated with the three-layer network. This experiment gives us insights for our follow-up work \cite{ADNN} on an adaptive network enhancement method which will study the problem of generating multi-layer networks, in terms of both width and depth, in order to approximate functions/PDEs of different characteristics accurately and efficiently. 

\section{Discussion and Conclusion}

%This paper studies a fundamental problem of approximating functions using two-layer ReLU NNs. Universal approximation theory guarantees that there exists a two-layer ReLU NN such that its approximation accuracy is within a given tolerance. 

This paper studies a fundamental question in machine learning on how to design the architecture of two-layer neural networks in order to approximate functions accurately and efficiently.
For a given function, we introduce and test an adaptive network enhancement (ANE) method that adaptively constructs a two-layer NN with a relatively small %nearly minimum 
number of neurons and parameters such that its approximation accuracy is within the prescribed tolerance. 
One of key components of the ANE method for the best least-squares approximation to a given function is the enhancement strategy which determines how many new neurons to be added, when the current approximation is not within the given accuracy. To address this issue, a global and a local network enhancement strategies are introduced and tested. The efficacy of the local enhancement strategy is demonstrated numerically for several test problems in this paper. Due to uncertainty of non-convex optimization, numerical results also show that the local strategy is better than the global one.  Nevertheless, efficiency and robustness of both the global and local enhancement strategies need further numerical and theoretical studies.

To disentangle the numerical integration error and network approximation error, an AMR method is proposed for automatically generating an integration mesh which adapts itself to improve the numerical integration accuracy. The AMR method presented in the paper is suitable for low dimensional problems and may be replaced by any adaptive integration procedure such as adaptive version of Monte Carlo, quasi-Monte Carlo, or sparse grid, etc. if a high dimensional problem is considered. Nevertheless, for a given function, how to adaptively choose a proper numerical integration in the context of NN functions remains open and requires further investigation.

Determining the values of the parameters of NNs is a problem in non-convex optimization which is computationally intensive and complicated and is a bottleneck in using NNs. Commonly used iterative solvers for optimization in NN applications are iterative methods of the gradient descent type. It is a common sense that it is extremely difficult, if not impossible, to develop a computationally feasible iterative solver that would converge to the desired global optimizer. This, in turn, implies the prominent importance of a close enough first approximation for all iterative solvers, as experienced in our numerical experiments. The method of continuation \cite{AlGe:90} is a common way to obtain a good initial and the ANE is a natural continuation process by itself with respect to the number of neurons.
In particular, weights and bias of newly added neurons are initialized based on the implicit physical partition of the domain $\Omega$ for the NN approximation at the previous network. This deterministic initialization strategy %, compared with other popular randomized initializers,  
ensures that the starting point of each iteration is always superior to the previous iteration when the network is enhanced, and plays an essential role in training the current network.

%This paper constructs such a network adaptively to meet the given approximation target. The proposed ANE method starts from a small network model, and it adaptively enhances its approximation power by adding neurons with proper initialization, which considers the implicit physical partition of a two-layer ReLU NN model and the corresponding local error indicators.  

Experimental results for functions exhibiting intersecting interface singularities or sharp interior layer like discontinuities show the efficacy of the propose method. In the second part of the paper \cite{LiuCai2020}, we extend the application of the proposed ANE method to elliptic partial differential equation with an underlying minimization principle.
%Moreover, they show a great promise of extending the proposed ANE method to other applications, e.g., solving for computationally challenging partial differential equations.
 
\bigskip
\bibliographystyle{elsarticle-num}
\bibliography{Reference}

\begin{thebibliography}{10}
\expandafter\ifx\csname url\endcsname\relax
  \def\url#1{\texttt{#1}}\fi
\expandafter\ifx\csname urlprefix\endcsname\relax\def\urlprefix{URL }\fi
\expandafter\ifx\csname href\endcsname\relax
  \def\href#1#2{#2} \def\path#1{#1}\fi

\bibitem{Hebb1949}
D.~O. Hebb, The organization of behavior: {A} neuropsychological theory, Wiley,
  New York, 1949.

\bibitem{Rosenblatt1958}
F.~Rosenblatt, The perceptron: A probabilistic model for information storage
  and organization in the brain, Psychological Review 65~(6) (1958) 386--408.

\bibitem{Cybenko1989}
G.~Cybenko, Approximation by superpositions of a sigmoidal function,
  Mathematics of Control, Signals, and Systems (MCSS) 2 (1989) 303--314.

\bibitem{HornikS1989}
K.~Hornik, M.~Stinchcombe, H.~White, Multilayer feedforward networks are
  universal approximators, Neural Networks 2 (1989) 359--366.

\bibitem{Petrushev1998}
P.~P. Petrushev, Approximation by ridge functions and neural networks, Siam
  Journal on Mathematical Analysis 30 (1998) 155--189.

\bibitem{pinkus1999}
A.~Pinkus, Approximation theory of the mlp model in nueral networks, Acta
  Numerica 8 (1999) 143--195.

\bibitem{LiuCai2020}
M.~Liu, Z.~Cai, Adaptive two-layer {R}e{LU} neural network {II}: Ritz
  approximation to elliptic {PDEs}, arXiv:2107.06459 [math.NA], (2021).

\bibitem{DickKuoSloan2014}
J.~Dick, F.~Kuo, I.~Sloan, High-dimensional integration - the quasi-monte carlo
  way, Acta Numerica 15 (2014) 133--288.

\bibitem{BungartzGriebel2004}
H.~J. Bungartz, M.~Griebel, Sparse grids, Acta Numerica 13 (2004) 1--123.

\bibitem{BoCuNo2018}
L.~Bottou, F.~E. Curtis, J.~Nocedal, Optimization methods for large-scale
  machine learning, SIAM Review 60 (2018) 223--311.

\bibitem{Schumaker}
L.~Schumaker, Spline Functions: Basic Theory, 1981.

\bibitem{jupp1978}
D.~Jupp, Approximation to data by splines with free knots, SIAM Journal on
  Numerical Analysis 15~(6) (1978) 328--343.

\bibitem{Baker85}
A.~J. Baker, On optimization aspects of a cfd finite element penalty algorithm,
  In: The Mathematics of Finite Element and Applications V (J. R. Whiteman,
  Ed.) (1985) 391--414.

\bibitem{Rice1969}
J.~R. Rice, The Approximation of Functions, Vol. 2, MA: Addison-Wesley, 1969.

\bibitem{Powell1968}
M.~Powell, On best l 2 spline approximations, Numerische Mathematik
  Differentialgleichungen Approximationstheorie (1968) 317--339.

\bibitem{chui1977}
C.~K. Chui, P.~W. Smith, J.~D. Ward, On the smoothness of best {L}2
  approximants from nonlinear spline manifolds, Math. Comput. 31 (1977) 17--23.

\bibitem{Dissanayake94}
M.~Dissanayake, N.~Phan-Thien, Neural network based approximations for solving
  partial differential equations, Communications in Numerical Methods in
  Engineering 10~(3) (1994) 195--201.

\bibitem{Sirignano18}
J.~Sirignano, K.~Spiliopoulos, {DGM}: A deep learning algorithm for solving
  partial differential equations, Journal of Computational Physics 375 (2018)
  1139--1364.

\bibitem{Karniadakis19}
M.~Raissia, P.~Perdikarisb, G.~Karniadakisa, Physics-informed neural networks:
  A deep learning framework for solving forward and inve, Journal of
  Computational Physics 378 (2019) 686–707.

\bibitem{cai2020}
Z.~Cai, J.~Chen, M.~Liu, X.~Liu, Deep least-squares methods: An unsupervised
  learning-based numerical method for solving elliptic pdes, Journal of
  Computational Physics 420 (2020) 109707.

\bibitem{Temlyakov2018}
V.~N. Temlyakov, The marcinkiewicz-type discretization theorems, Constructive
  Approximation 48 (2018) 337--369.

\bibitem{PCA1901}
K.~Pearson, On lines and planes of closest fit to systems of points in space,
  The London, Edinburgh, and Dublin Philosophical Magazine and Journal of
  Science 2~(11) (1901) 559--572.

\bibitem{kingma2014adam}
D.~P. Kingma, J.~Ba, Adam: A method for stochastic optimization, arXiv preprint
  arXiv:1412.6980 (2014).

\bibitem{he2018relu}
J.~He, L.~Li, J.~Xu, C.~Zheng, Relu deep neural networks and linear finite
  elements, Journal of Computational Mathematics 38~(3) (2020) 502--527.

\bibitem{Morin2002}
P.~Morin, R.~H. Nochetto, K.~G. Siebert, Convergence of adaptive finite element
  methods, SIAM Review 44~(4) (2002) 631--658.

\bibitem{Cai2009}
Z.~Cai, S.~Zhang, Recovery-based error estimator for interface problems:
  Conforming linear elements, SIAM Journal on Numerical Analysis 47~(3) (2009)
  2132–2156.

\bibitem{ADNN}
Z.~Cai, J.~Chen, M.~Liu, Self-adaptive deep neural network: numerical
  approximation to functions and {PDEs}, arXiv:2109.02839 [math.NA], (2021).

\bibitem{AlGe:90}
E.~Allgower, K.~Georg, Numerical Continuation Methods, Springer-Verlag, Berlin
  and Heidelberg, 1990.

\end{thebibliography}

\end{document}